\documentclass[12pt]{amsart}

\usepackage{hyperref}

\usepackage{soul}

\def\cal{\relax}

\usepackage{DLdef1}
\let\ssec\subsection

\renewcommand {\ssbegin}[1]
 {\refstepcounter{subsection}
 \def \secno {\gdef \secno {}{\ssecfont
\thesubsection.\hskip 2ex}%
 }%
 \begin{#1}}

 \def\mmat #1,#2,#3,#4,{\text{\small\arraycolsep=3pt $
\begin{pmatrix}#1&#2\\#3&#4\end{pmatrix}$}}

\newcommand{\un}{\underline{N}}

\newcommand{\del}{\partial}

\begin{document}


\thispagestyle{empty}

\title[Invariant differential operators]{Invariant differential operators\\
in positive characteristic}

\author{Sofiane Bouarroudj${}^{1*}$, Dimitry Leites${}^2$}

\address{${}^1$New York University Abu Dhabi,
Division of Science and Mathematics, P.O. Box 129188, United Arab
Emirates; sofiane.bouarroudj@nyu.edu (${}^*$ Corresponding author) \\
${}^2$Department of Mathematics, Stockholm University, SE-106 91
Stockholm, Sweden; mleites@math.su.se}



\begin{abstract} We consider an analog of the problem Veblen formulated in 1928 at
the IMC: classify invariant differential operators between ``natural
objects" (spaces of either tensor fields, or jets, in modern terms)
over a real manifold of any dimension. For unary operators, the
problem was solved by Rudakov (no nonscalar operators except the
exterior differential); for binary ones, by Grozman (there are no
operators of orders higher than 3, operators of order 2 and 3 are,
bar an exception in dimension 1, compositions of order 1 operators
which, up to dualization and permutation of arguments, form 8
families). In dimension one, Grozman discovered an indecomposable
selfdual operator of order 3 that does not exist in higher
dimensions. We solve Veblen's problem in the 1-dimensional case over
the ground field of positive characteristic. In addition to analogs
of the Berezin integral (strangely overlooked so far) and binary
operators constructed from them, we discovered two more (up to
dualization) types of indecomposable operators of however high
order: analogs of the Grozman operator and a completely new type of
operators.

\end{abstract}




\maketitle

\section{Introduction}

In the representation theory of a given algebra or group, the usual
ultimate goal is a clear description of its irreducible modules to
begin with, indecomposable ones at the next step. Other goals are
sometimes also natural and reachable, cf. those we discuss here and
in review \cite{BF}.

\ssec{Discussion of open problems} The history of mathematics shows
that problems with a nice (at least in some sense, for example,
short) answer often turn out to be ``reasonable". I.~Gelfand used to
say that, the other way round, easy-to-formulate lists are answers
to ``reasonable" problems  and advised to try to formulate such
problems once we got a short answer.

\textbf{Hereafter, $\Kee$ is any field of characteristic $p>0$,
unless otherwise specified}.

Speaking of Lie algebras or superalgebras, and their representations
over $\Kee$, which of them is it natural to consider? In \cite{LL},
Deligne suggests that one should begin with restricted ones: unlike
nonrestricted ones, the restricted Lie (super)algebras correspond to
groups, in other words: to geometry. Certain problems concerning
nonrestricted algebras are also natural (although tough), and have a
short answer, e.g., classification of simple Lie algebras for $p>3$
(for its long proof, see \cite{S,BGP}).

Rudakov and Shafarevich \cite{RSh} were the first to describe ALL,
not only restricted, irreducible representations of  $\fsl(2)$ for
$p>2$. Dolotkazin solved the same problem for $p=2$, see \cite{Do};
more precisely, he described the irreducible representations of the
simple 3-dimensional Lie algebra $\fo^{(1)}(3)$. The difference of
Dolotkazin's problem from that considered in \cite{RSh} is that,
unlike $\fsl(2)$ for $p>2$, the Lie algebra $\fo^{(1)}(3)$ is not
restricted. These two results show that the description of all
irreducible representations looks feasible, to an extent, at least
for Lie (super)algebras with indecomposable Cartan matrix or their
simple ``relatives" (for the classification of both types of
(super)algebras over $\Kee$, see \cite{BGL}).

The problem we are solving here for $p>0$ resembles its super
version \textbf{over $\Cee$}; we recall the known results. Simple
Lie (super)algebras of vector fields with polynomial or formal
coefficients constitute another class of interesting Lie
(super)algebras. \textbf{Over $\Cee$}, they are classified and their
continuous (in the $x$-adic topology) irreducible representations
with a vacuum (highest or lowest) weight vector are classified as
well in the following cases:

$\bullet$ For the four series of simple Lie algebras, and their
direct super versions, the answer is as follows, see the review
\cite{GLS}: the module $T(V)$ of tensor fields coinduced from any
irreducible module $V$ with lowest weight vector over the subalgebra
of linear vector fields is irreducible, unless $T(V)\simeq\Omega^k$,
the module of exterior $k$-forms.

Rudakov proved that the exterior differential is the only
(nonscalar) operator between the modules of exterior differential
forms invariant with respect to all changes of indeterminates, or,
equivalently, for the general vectorial Lie algebra. Together with
the Poincar\'e lemma (stating that the sequence \eqref{ext} is exact
over any star-shaped domain) this yields a complete description of
continuous irreducible modules with lowest weight.

$\bullet$ For the remaining 3 series of simple vectorial Lie
algebras (and analogous simple Lie superalgebras \textbf{in their
standard grading}) the answer is essentially the same.

$\bullet$ For simple vectorial Lie superalgebras in their
\textbf{nonstandard} grading (listed in \cite{LS}), and for the
modules of tensor fields with \textbf{infinite-dimensional fibers}
(i.e.,  $\dim V=\infty$) even over finite-dimensional simple
vectorial Lie superalgebras, many never previously known
differential operators appear, see \cite{GLS,Lint}.

$\bullet$ In the purely odd case the Berezin integral has to be
added to the list of invariant differential operators, see
eq.~\eqref{intDR}.

$\bullet$ \textbf{Over $\Kee$}, it is tempting to conjecture that,
at least for $p$ ``sufficiently large", the description of
irreducible modules over vectorial Lie algebras (at least
$\Zee$-graded, and ``close to simple" ones) the description of
irreducible modules would resemble that given by Rudakov over
$\Cee$. For irreducible  restricted modules over restricted simple
vectorial Lie algebras of the 4 classical series and $p\geq 3$, this
is almost so, i.e., the spaces of differential $k$-forms areg the
only reducible spaces of tensor fields, as proven by Krylyuk, see
\cite{Kry}, although the complex of exterior differential forms is
not exact, see Theorem~\ref{Kryl}.

There are, however, natural filtered deforms of simple nonrestricted
vectorial Lie algebras; for $p\leq 5$, there are new types of simple
vectorial Lie algebras and no \textbf{description of their
irreducible modules is known; in particular, no description of
intertwining operators between modules of tensor fields --- the most
natural ones over these Lie algebras}. Description of this type has
a history.

\sssec{History: Veblen's problem} In 1928, at the IMC, O.~Veblen
formulated a problem which was later reformulated in more
comprehensible terms by A.~Kirillov and further reformulated as a
purely algebraic problem by J.~Bernstein who interpreted Rudakov's
solution of Veblen's problem for unary operators; for setting in
modern words and review, in particular, for a superization of
Veblen's problem, see \cite{GLS}.

Let $\fvect(m)$ be the Lie algebra of polynomial vector fields (over
a ground field of characteristic 0, say $\Cee$). Assuming the notion
of the $\fvect(m)$-module $T(V)$ of (formal) tensor fields of type
$V$, where $V$ is a $\fgl(m)$-module with lowest weight vector, is
known, see \cite{GLS}. Let us briefly recall the results concerning
the \textbf{nonscalar} unary and binary invariant differential
operators between spaces of tensor fields, although we only need the
simplest version of spaces $T(V)$, namely, the spaces of weighted
densities $\cF_a:=T(a\, \tr)$, where $\tr$ is the 1-dimensional
$\fgl(m)$-module given by the trace (supertrace for $\fgl(m|k)$) and
$a$ is in the ground field. We assume that modules $T(V)$ in what
follows are coinduced from \textbf{irreducible} modules $V$; this is
certainly so for modules $\cF_a$. (In Veblen's problem for simple
vectorial Lie algebras over $\Cee$, it suffices to consider only
irreducible modules $V$; for its superization and over $\Kee$, one
has to consider indecomposable modules $V$, but one has to begin
with irreducible ones, anyway.)

\paragraph{Unary operators} The only unary $\fvect(m)$-invariant
differential operators $T(V)\tto T(W)$ are the exterior differential
operators (so, practically, there is just one such operator) in the
de Rham complex, where $\Omega^i=T(E^i(\id))$, $E^i$ is the functor
of raising to the $i$th exterior power and $\id$ is the tautological
$m$-dimensional $\fgl(m)$-module:
\begin{equation}
\label{ext}
0\tto\Cee\tto\Omega^0\stackrel{d}{\tto}\Omega^1\stackrel{d}{\tto}\dots\stackrel{d}{\tto}
\Omega^{m-1} \stackrel{d}{\tto}\Omega^m\tto 0.
\end{equation}

In super setting, where $\sdim(\id)=m|k$, the analog of the complex
\eqref{ext} is infinite without the ``top term"
\begin{equation}
\label{extSup}
0\tto\Cee\tto\Omega^0\stackrel{d}{\tto}\Omega^1\stackrel{d}{\tto}\Omega^2
\stackrel{d}{\tto}\dots.
\end{equation}
its dualization brings the complex of \textit{integrable} forms
which \textbf{in the purely odd case of superdimension $0|k$}
terminates with an order-$k$ \textbf{differential} operator called
the \textit{Berezin integral}:
\begin{equation}
\label{intDR}
\renewcommand{\arraystretch}{1.4}
\begin{array}{ll}
\dots\stackrel{d}{\tto}\Sigma_{-1}(m|k)\stackrel{d}{\tto}
\Sigma_{0}(m|k)\stackrel{d}{\tto}0&\text{for $m\neq 0$},\\
\dots\stackrel{d}{\tto}\Sigma_{-1}
(0|k)\stackrel{d}{\tto}\Sigma_{0}(0|k)
 \stackrel{\int}{\tto}\Kee\tto 0&\text{for $m= 0$}.\\ \end{array}
\end{equation}
Over $\Cee$, the de Rham complex \eqref{extSup} is exact, while the
complex of integrable forms \eqref{intDR} has 1-dimensional
cohomology in the $(-m)$th  term; the cohomology is spanned, see
\cite{BL,MaGF,Del}, by (here the $u_i$ are the even indeterminates,
the $\theta_j$ are the odd ones for clarity; in Notation~\ref{Nota}
all indeterminates are called $u$)
\begin{equation}
\label{intCoh}\theta_1\dots\theta_k\del_{u_1}\dots
\del_{u_m}\vvol(u|\theta)
\end{equation}
and the super version of the volume element
$\vvol(u)=du_1\wedge\dots\wedge du_m$ is  the class
$\vvol(u|\theta)$ of the element
\begin{equation}
\label{vvol} du_1\wedge\dots\wedge du_m\otimes
\del_{\theta_1}\wedge\dots\wedge\del_{\theta_k}
\end{equation}
in the indecomposable $\fgl(V)=\fgl(m|k)$-module generated by the
element \eqref{vvol} modulo the codimension 1 submodule, see
\cite{Del,LSoS}.

\paragraph{Binary operators} These are classified by P.~Grozman, see
\cite{GInv}. Recall that the formal dual of $T(V)$ on the
$m$-dimensional space is not $T(V^*)$ but $T(V^*\otimes \tr)$.
Hence, with every unary operator $D\colon T(V)\tto T(W)$ there
exists its dual $D^*\colon T(W^*\otimes \tr)\tto T(V^*\otimes \tr)$.

In particular, on the line,  with every operator $D\colon \cF_a\tto
\cF_{a+\deg D}$ of order $\deg D$, there exists its dual $D^*\colon
\cF_{1-a-\deg D}\tto \cF_{1-a}$ of the same order.

Similarly, with every binary operator $D\colon \cF_a\otimes
\cF_b\tto \cF_{a+b+\deg D}$, there exist the two its duals
$D^{*1}\colon \cF_{1-a-b-\deg D}\otimes \cF_b\tto \cF_{1-a}$ and
$D^{*2}\colon \cF_a\otimes \cF_{1-a-b-\deg D}\tto \cF_{1-b}$, as
well as the operator $r(D)\colon \cF_b\otimes \cF_a\tto
\cF_{a+b+\deg D}$, given by the interchange of the arguments
\[
r(D)(X,Y)\colon =D(Y,X)\text{~~for any $X\in \cF_b$ and
$Y\in\cF_a$.}
\]

\textbf{Summary}: Up to dualizations and interchange of the
arguments, there are 8 series of invariant operators of order 1, all
invariant differential operators of order 2 and 3 are compositions
of order-1 operators, except for the case $m=1$, where there is an
indecomposable operator $Gz$ of order $3$ discovered by Grozman; for
explicit expressions, see \cite{GInv} and Subsec.~\ref{BiCee}. There
are no operators of order $>3$.

\textbf{Off-topic}: 1) Among all 8 series of invariant operators of
order 1, there are 4 series that determine an associative algebra
and 3 series determine a Lie (super)algebra structure on their
domain of definition. These Lie (super)algebras are either simple or
simple modulo center or contain a simple algebra of codimension 1:
the corresponding products are the Poisson bracket, the Schouten
bracket (more popular now under the name anti-bracket) and a
deformed Schouten bracket, see \cite{GInv, LSh}.

2) Kirillov gave an example of the symbol of an invariant with
respect to all coordinate changes \textit{nonlocal} bilinear
operator on the circle but did not identify the operator itself; in
\cite{IoMo} this operator is described in terms of the
integral=residue; it is the only nonlocal invariant bilinear
operator on the circle. In characteristic $p>0$, same as on
superpoints, analogs of various operators nonlocal over $\Ree$
become local because if the space of functions is
finite-dimensional, then any operator in this space can be viewed as
differential, i.e., composition of derivations and multiplications
by some functions.

\paragraph{Invariant differential operators of higher arity}
The only $k$-ary invariant operators classified so far are
antisymmetric ones on $\Cee$ , see \cite{FeFu}, and ternary
operators between the spaces of weighted densities on $\Cee^n$, see
\cite{Bj}.

\ssec{Main result} The classification (Theorem~\ref{Th}) of binary
$\fvect(1;\un)$-invariant operators between modules of weighted
densities; in particular, we discovered three series of
indecomposable operators of however high order: 1) related to the
analog of Berezin integral, 2) an analog of the Grozman operator,
and 3) a completely new series of operators, $Bj$.

\ssec{Open problems} 1) We conjecture that the Grozman operator has
something to do with the contact structure since $\fvect(1)=\fk(1)$.
Over $\Cee$, it seems feasible to classify binary
$\fk(2n+1)$-invariant differential operators for $n>0$ and complete
the partial result for binary $\fh(2n)$-invariant differential
operators, solved only for $n=1$, see \cite{Ghm}.

2) The only feasible super version of the problem we consider here
is, it seems, classification of the $\fk(1;\un|k)$-invariant
operators for $k=1$ or 2 for which the coinduced modules with lowest
weight vectors are of the form $\cF_a$. For the other pairs $m,k$
the analog of Veblen's problem does not seem to be feasible.

3) Prove Conjecture~\ref{Conj}.

\section{Veblen's problem over fields of
characteristic $p>0$}

\ssec{Notation}\label{Nota} For an $m+k$-tuple of nonnegative
integers $\underline{r}=(r_1, \ldots , r_a)$, where $r_i=0$ or 1 for
$i>m$, we introduce parity by setting $p(u_i)=\ev$ for $i\leq m$ and
$p(u_i)=\od$ for $i> m$. Set further
\begin{equation*}
\label{1.} u_i^{(r_{i})} := \frac{x_i^{r_{i}}}{r_i!}\quad
\text{and}\quad u^{(\underline{r})} := \prod\limits_{1\leq i\leq a}
u_i^{(r_{i})}.
\end{equation*}
The idea is to formally replace fractions with $r_i!$ in
denominators by inseparable symbols $u_i^{(r_{i})}$ which are
well-defined over fields of characteristic $p$ because the structure
constants are integers one can consider modulo $p$, in particular,
if $k=0$, the factor in the first parentheses below is equal to 1:
\begin{equation*}
\label{divp}
\renewcommand{\arraystretch}{1.4}
\begin{array}{l} u^{(\underline{r})} \cdot u^{(\underline{s})} =
\left(\mathop{\prod}\limits_{m+1\leq i\leq n}
\min(1,2-r_i-s_i)\cdot(-1)^{\mathop{\sum}\limits_{m<i<j\leq a}
r_js_i}\right) \cdot \binom {\underline{r} + \underline{s}}
{\underline{r}}
u^{(\underline{r} + \underline{s})}, \\
\text{where}\quad \binom {\underline{r} + \underline{s}}
{\underline{r}}:=\mathop{\prod}\limits_{1\leq i\leq m}\binom {r_{i}
+ s_{i}} {r_{i}}.
\end{array}
\end{equation*}
For an $m+k$-tuple of positive integers $\underline{N} = (N_1,...,
N_m, 1, \dots, 1)$, set $\underline{N}_{ev} = (N_1,..., N_m)$. The
following supercommutative superalgebra is the analog of the algebra
of functions when $p>0$:
\begin{equation*}
\label{u;N}
\renewcommand{\arraystretch}{1.4}
\begin{array}{l}
\cO(m;\un_{ev}|k):=\Kee[u;\un]:=
\Span_{\Kee}(u^{(\underline{r})}\mid r_i < p^{N_{i}}\;\text{ for
$i\leq m$, $0\leq r_i \leq 1$ for $i> m$}).\end{array}
\end{equation*}
The Lie superalgebra $\fvect(m;\un|k)$ of its \textit{distinguished}
derivations consists of the vector fields of the form $\sum
f_i(u)\del_i$, where ${f_i\in\cO(m;\un_{ev}|k)}$ and $\del_i$ are
\textit{distinguished} partial derivatives defined by the condition
${\del_j(u_i^{(r_{i})})=\delta_{ji}u_i^{(r_{i}-1)}}$.

\ssec{Unary operators: De Rham complex for $k=0$ (no super)}
\sssbegin{Theorem}[\cite{Kry}]\label{Kryl} Set
\begin{equation*}\label{tau}
\tau(\underline{N})=(p^{N_1}-1, \dots, p^{N_m}-1).
\end{equation*}
Denote: $Z^i(m;\underline{N}):=\{\omega
\in\Omega^i(m;\underline{N})\mid d\omega=0\}$ and
$B^i(m;\underline{N}):=\{d\omega \mid \omega
\in\Omega^{i-1}(m;\underline{N})\}$.

On the $m$-dimensional space over $\Kee$, the sequence
\begin{equation*}
\label{ext-p}
\renewcommand{\arraystretch}{1.4}
\begin{array}{l}
0\tto\Kee\tto\Omega^0(\underline{m;N})\stackrel{d}{\tto}
\Omega^1(m;\underline{N})\stackrel{d}{\tto}\Omega^2(m;\underline{N})
\stackrel{d}{\tto}\dots \stackrel{d}{\tto}\Omega^m(m;\underline{N})
\stackrel{d}{\tto}0 \end{array}
\end{equation*}
is not exact: the space
$H^a(m;\underline{N}):=Z^a(m;\underline{N})/B^a(m;\underline{N})$ is
spanned by the elements
\begin{equation*}
\label{coh} u_{i_1}^{(\tau(\underline{N})_{i_1})}\dots
u_{i_k}^{(\tau(\underline{N})_{i_k})}du_{i_1}\dots du_{i_k}\text{
where $a=i_1+\dots +i_k$}.
\end{equation*}
In other words, the supercommutative superalgebra
$H^{\bcdot}(m;\underline{N})$ is generated by
$H^1(m;\underline{N})$.
\end{Theorem}

\begin{proof}Induction on $m$. For $m=1$, this is obvious. \end{proof}

\sssbegin{Corollary}\label{Cor} There is an invariant differential
operator of order $|\tau(\un)|:=(\sum p^{\un_i})-m$:
\begin{equation*}
\label{int-p}
\renewcommand{\arraystretch}{1.4}
\begin{array}{l}
\Omega^m(m;\underline{N}) \stackrel{\int}{\tto}\Kee; \end{array}
\end{equation*}
this $\int$ is a, depending on $\un$, analog of the Berezin integral
defined as follows:
\[
\int  f(u)\vvol(u)=\text{the coefficient of
$u^{\tau(\un)}=u_{1}^{(\tau(\underline{N})_{1})}\dots
u_{m}^{(\tau(\underline{N})_{m})}$, where
$\vvol(u)=du_1\wedge\dots\wedge du_m$}.
\]
\end{Corollary}

\parbegin{Remark} We did not find the statement of the above
Corollary in the literature; so formulation of the fact that there
exists an integral depending on $\un$ seems to be a new
result.\end{Remark}

\ssec{Unary operators: De Rham complex for $k\neq 0$, i.e., in super
setting} On the $m|k$-dimensional space over $\Kee$,  the de Rham
complex goes ad infinitum
\begin{equation*}
\label{superDR}
\renewcommand{\arraystretch}{1.4}
\begin{array}{l}
0\tto\Kee\tto\Omega^0(m;\un|k)\stackrel{d}{\tto}
\Omega^1(m;\un|k)\stackrel{d}{\tto}\Omega^2(m;\un|k)
\stackrel{d}{\tto}\dots \end{array}
\end{equation*}

The dual to \eqref{superDR} complex of \textit{integrable} (Deligne
calls them \textit{integral}) forms is as follows, where
$\Sigma_{0}:=\Omega^0\vvol(u|\theta)$:
\begin{equation}
\label{intDRp>0}
\renewcommand{\arraystretch}{1.4}
\begin{array}{l}
\dots\stackrel{d}{\tto}\Sigma_{-1}
(m;\un|k)\stackrel{d}{\tto}\Sigma_{0}(m;\un|k)
 \stackrel{\int}{\tto}\Kee\tto 0. \end{array}
\end{equation}

Observe that, for $p>0$, the integral is defined for $m\neq0$ as
well, whereas for $p=0$, there is no such an invariant
\textit{differential} operator.

\sssbegin{Conjecture}\label{Conj} The only
$\fvect(m;\un|k)$-invariant unary operators between modules of
tensor fields coinduced from irreducible $\fgl(m|k)$-modules  are
the exterior operator $d$ and the Berezin integral $\int$.
\end{Conjecture}


\ssec{The binary differential operators over $\Cee$ for $m=1$ (recapitulation,
see \cite{GInv})}\label{BiCee} For $m=1$, we have
$\cF_a:=\{f(u)(du)^a\mid f(u)\in \Omega^0=\cF_0\}$ on which
the Lie derivative acts as follows:
\begin{equation}\label{F_a}
L_D(f(u)(du)^a)=(D(f)+afg')(du)^a)\text{~~for any
$D=g(u)\partial_u\in\fvect(1)$}.
\end{equation}

For any bilinear differential
operator
\[
D\colon \cF_a\otimes \cF_b\tto \cF_{a+b+\deg D},\ \ \
f(du)^a\otimes g(du)^b\longmapsto D(f,g)(du)^{a+b+\deg D},
\] it suffices to
indicate $(a,b)$ and $D(f,g)$. The operators are listed up to
proportionality.

\sssec{Order 1 operators} For generic values $(a,b)$, the invariant
operators form a 1-dimensional space. For $(a,b)=(0,0)$, and only in
this case, we have a 2-dimensional space of operators.

\begin{equation}
\label{.1} \renewcommand{\arraystretch}{1.4}
\begin{tabular}{|l|l|l|}
\hline
$(a,b)$&$D(f,g)$&\text{comment} \\
\hline
$(0,0)$&$P_{00}\colon \alpha f'g+\beta fg'$&\text{for any $\alpha,\beta\in\Cee$}\\
\hline $(a,b)$&$\{\cdot,\cdot\}_{P.B.}\colon =afg'-bf'g$
&\text{Poisson bracket in two
even indeterminates $u$ and $du$}\\
$(-1,-1)$&$\{\cdot,\cdot\}\colon =fg'-f'g$ &\text{contact bracket of
generating
functions $f\longmapsto K_f:=f\nfrac{d}{du}$}\\
\hline\end{tabular}
\end{equation}

\sssec{Order 2 operators}
\begin{equation}
\label{.2} \begin{tabular}{|l|l|l|}
\hline$(a,b)$&$D(f,g)$&\text{comment}\\
\hline $(0,b)$
&$f'g'-bf''g$&$\{df,g(du)^b\}_{P.B.}$\\
$(a,0)$&$afg''-f'g'$&$\{f(du)^a,dg\}_{P.B.}$ \\
$(a,-1-a)$&$afg''+(2a+1)f'g'+(a+1)f''g$&$(\{f(du)^a,g(du)^{-1-a}\}_{P.B.})'$ \\
\hline\end{tabular}
\end{equation}

\sssec{Order 3 operators}
\begin{equation}
\label{.3} \begin{tabular}{|l|l|}
\hline$(a,b)$&$D(f,g)$\\
\hline
$(0,0)$&$T_1(f,g)=\{f',g'\}$\\
$(-2,0)$&$T_1^{*1}(f,g)=f'g''+3f''g'+ 2f'''g$\\
$(0,-2)$&$T_1^{*2}(f,g)=f''g'+3f'g''+ 2fg'''$\\
\hline

$\left(-\nfrac23,-\nfrac23\right)$&$
Gz(f,g)=2\det\begin{pmatrix}f&g\\
f'''&g'''\end{pmatrix}+3\det\begin{pmatrix}f'&g'\\
f''&g''\end{pmatrix}$\\
\hline\end{tabular}
\end{equation}

\ssec{New binary differential operators: $p>0$ and $m=1$} For a
fixed $\un=N$, we list $\fvect(1;\un)$-invariant operators up to
interchange of arguments.

Recall that $\int\colon  \Omega^1=\cF_1\tto\Kee\subset\cF_0$ is the
order $p^N-1$ operator, see Corollary~$\ref{Cor}$ and $1-p^N\equiv
1\pmod p$. Therefore, for any $p>0$ and $\un$, there are
$\fvect(1;\un)$-invariant bilinear operators acting in the spaces
$\cF_a\otimes\cF_b$ in the following cases:
\begin{equation}\label{ind1}
\begin{tabular}{|l|c|c|}
\hline
operator$(a,b)$&order&case\\
\hline $\int\otimes \id(1,a)$,\ \ $(\int\otimes \id)^{*1}(1-a,a)$,\
\ $(\int\otimes
\id)^{*2}(1,1-a)$&$k=p^N-1$&1\\

$\int\otimes d=(\int\otimes d)^{*2}(1, 0)$,\ \ $(\int\otimes d)^{*1}(0,0)$,\ \
$r(\int\otimes d)(0,1)$ &$k=p^N$&2\\

$\int\otimes \int=(\int\otimes \int)^{*1}=(\int\otimes
\int)^{*2}(1,1)$&$k=2(p^N-1)$&3\\
\hline
\end{tabular}
\end{equation}
The operator $(\int\otimes \id)^{*1}(1-a,a)$ is of the form ($L$ for
``long" expression)
\begin{equation}\label{LongOp}
L_{p,N}(f,g)=\sum_{i=0}^{\lfloor (k-1)/2\rfloor} (-1)^{i}
\left ( f^{(i)}g^{(k-i)}+ f^{(k-i)}g^{(i)}\right ) +
\begin{cases} (-1)^{k/2}f^{(k/2)}g^{(k/2)}&\text{if $p>2$}\\
&\text{if $p=2$}.\end{cases}
\end{equation}

\sssbegin{Lemma}\label{leQ} For  $(a,b)=(1,1)$, there exists a
$\fvect(1;\un)$-invariant bilinear operator $Bj_{p,m}$ of order
$k=p^m-1$ for any $m\leq N$ given by  the expression
\begin{equation}\label{ind2}
Bj_{p,m}(f,g)=\det
\begin{pmatrix}f &g\\
f^{(k)}&g^{(k)}\end{pmatrix}.
\end{equation}

Its duals for $(a,b)=(1,0)$ and $(0,1)$, respectively, are as
follows:

\begin{equation}\label{ind3}
Bj_{p,m}^{*2}(f,g)=f g^{(k)}+\sum_{i=1}^{\lfloor (k-1)/2\rfloor}
(-1)^{i} \left ( f^{(i)}g^{(k-i)}+ f^{(k-i)}g^{(i)}\right ) +
\begin{cases} (-1)^{k/2}f^{(k/2)}g^{(k/2)}&\text{if $p>2$}\\
&\text{if $p=2$ }\end{cases}
\end{equation}
and
\begin{equation}\label{ind4}
Bj_{p,m}^{*1}(f,g)=f^{(k)} g+\sum_{i=1}^{\lfloor (k-1)/2\rfloor}
(-1)^{i} \left ( f^{(i)}g^{(k-i)}+ f^{(k-i)}g^{(i)}\right ) +
\begin{cases} (-1)^{k/2}f^{(k/2)}g^{(k/2)}&\text{if $p>2$}\\
&\text{if $p=2$ }.\end{cases}
\end{equation}
\end{Lemma}

\begin{proof}Every order $k$ bilinear differential operator $D$ is of the form
\[
(f,g)\longmapsto\sum_{i+j=k} \alpha_{i,j} f^{(i)}
g^{(j)},\text{~~where $\alpha_{i,j}\in\Kee$}.
\]
For $N$ big enough (such that $p^N>\deg D+1$) and a given pair
$(a,b)$, the $\fvect(1;\un)$-invariance of $D$ implies the following
system of equations
\begin{equation*}
\label{InvEq1}
\left ( a \binom{i}{r-1} + \binom{i}{r}  \right )\alpha_{i, k-i}+
\left ( b \binom{k-i+r-1}{r-1} + \binom{k-i+r-1}{r}  \right
)\alpha_{i-r+1, k-i+r-1}=0,
\end{equation*}
for $i=1,\ldots, k$ and $2\leq r \leq k+1$. In the case of Lemma,
$(a,b)=(1,1)$; besides $\alpha_{i,j}=0$ whenever $i+j\not=k$. The
above system becomes
\begin{equation}
\label{InvEq2} \binom{i+1}{r} \alpha_{i, k-i}+  \binom{k-i+r}{r}
\alpha_{i-r+1, k-i+r-1}=0\text{~~for $i=1,\ldots, k$ and $2\leq r
\leq k+1$}.
\end{equation}
If $r=k+1$, then $i$ must be equal to $k$, so we get the condition
$\alpha_{0,k}+\alpha_{k,0}=0.$ If $r<k+1$, then the above equations
are either identically zero (since $\alpha_{i,j}=0$ whenever
$i+j\not=k$) or of the form
\[
\binom{k+1}{r} \alpha_{k, 0}=0 \text{ and } \binom{k+1}{r}
\alpha_{0, k}=0.
\]
But $\binom{k+1}{r}=0\, (\text{mod } p)$, since $k+1=p^N$.

The claim on dual operators is subject to a direct verification.
\end{proof}

\sssec{The transvectants (also known as Cohen-Rankin brackets)} Recall that for any $p\neq 2$, there is an
embedding $\fsl(2)\hookrightarrow \fvect(1;\un)$. (For $p=2$, the
analog of this embedding is the embedding
$\fo^{(1)}(3)\hookrightarrow \fvect(1;\un)$ for $\un>1$.) The
\textit{transvectants} are $\fsl(2)$-invariant bilinear differential
operators discovered by Gordan. As a matter of fact, there is a
misprint in the explicit expressions of the transvectants
$J_{a,b}^k$ in \cite{BjOv}; the corrected expression is (up to a
numerical factor)
\begin{equation}\label{ind5} J_{a,b}^k(f,g)= \sum_{0\leq i\leq k}(-1)^i
\binom{2a +k-1}{k-i}\binom{2b+k-1}{i} f^{(i)} g^{(k-i)}.
\end{equation}
This expression is for generic $a$ and $b$, and $p\not =2$. For the
exceptional values of $a$ and $b$ over $\Cee$, see \cite{Bj2},
Prop.2. In this paper, we are mainly interested in the explicit
expressions
of the transvectants for $p=2$. 

\sssbegin{Lemma}\label{CjGz} For  $(a,b)=(1,1)$, there exists a
$\fvect(1;\un)$-invariant bilinear operator of order $k=p^m-2$ for
any $m\leq N$ if $p\not =2$ whereas for $p=2$ it should be any
$m\leq N+1$ (because of the term $f^{(k/2)}g^{(k/2)}$) given by the
expression
\begin{equation}\label{ind6}
Gz_{p,m}(f,g)=\sum_{i=0}^{\lfloor{(k-1)/2\rfloor}} (-1)^{i}\det
\begin{pmatrix}f^{(i)}&g^{(i)}\\
f^{(k-i)}&g^{(k-i)}\end{pmatrix}+
\begin{cases}
f^{(k/2)}g^{(k/2)}&\text{if $p=2$}\\
&\text{otherwise}.\\
\end{cases}
\end{equation}
\end{Lemma}
This operator coincides with the Grozman operator if $p=5$ and
$m=1$. Observe that $Gz$ is selfdual over any field it is defined,
but $Gz\circ (d\otimes d)$ is not selfdual.

\begin{proof} Let, for simplicity, $p\not =2$ (for $p=2$ the proof is more or less the
same). Since $(a,b)=(1,1)$, the invariance with respect to
$\fvect(1;\un)$ is equivalent to the system (\ref{InvEq2}). For
$i=k$ or $i=0$, we have $\alpha_{k,0}+\alpha_{0,k}=0$ which is
certainly satisfied. Now, for $i\not=0,k$, the system becomes
\[
\binom{i+1}{r} (-1)^i+  \binom{k-i+r}{r}  (-1)^{i-r+1}=0
\text{~~for $i=1,\ldots, k$ and $2\leq r \leq k+1$}.
\]

To show that the system is satisfied, we proceed by induction. If
$i=1$, then $r$ must be equal to 2. The result follows since
\[
\begin{array}{lcl}
\binom{k+1}{2} & =& \binom{p^m-1}{2} = 1 \text{ mod } (p).\\[3mm]
\end{array}
\]
Now suppose that the equality is
true at $i$ for every $r$ such that $r\leq i+1$. It follows that
\[
\begin{array}{lcl}
\binom{k-(i+1)+r}{r} & =&\binom{k-i+r}{r}-\binom{k-i+r-1}{r-1} \text{ `Pascal's triangle' }\\[3mm]
& =&(-1)^r\binom{i+1}{r}-(-1)^{r-1}\binom{i+1}{r-1} \text{ mod } (p) \text{ induction hypothesis}\\[3mm]
& =&(-1)^r \left (\binom{i+1}{r}+\binom{i+1}{r-1}\right ) \text{ mod } (p) \\[3mm]
& =&(-1)^r \binom{i+2}{r} \text{ mod } (p) \text{ `Pascal's triangle' } \hfill \qed\\
\end{array}
\]

\noqed\end{proof}

\ssbegin{Theorem}\label{Th} Up to the interchange of arguments, the
indecomposable $\fvect(1;\un)$-invariant bilinear operators $D\colon
\cF_a\otimes\cF_b\tto\cF_c$ are only those of the form
$\eqref{ind1}$--$\eqref{ind6}$. For the complete list of
$\fvect(1;\un)$-invariant bilinear operators of order $\leq 7$, see
tables $\eqref{t1}$--$\eqref{t7}$.
\end{Theorem}

Although we have found all invariant operators of order $\leq 20$,
the complete list is not that interesting, we think, unlike the
indecomposable operators. We illustrate the answer with a part of
the complete list.

\begin{proof} Computer-aided study with the aid of
Grozman's \textit{Mathematica}-based code \textit{SuperLie}, see
\cite{Gr}:  (1) We fix $\deg D\in\{1, 2, \dots, 20\}$ and
$p\in\{2,3,5,\dots, 19\}$. (2) For this $\deg D$ and $p$, we look
for which $\un$ and $(a,b)$ there is a $\fvect(1;\un)$-invariant
operator.
\end{proof}


$\bullet$ \emph{Order 1 operators} are the same as for $p=0$, but
for $p=2$ we have
\begin{equation}
\label{t1} \renewcommand{\arraystretch}{1.4}
\begin{tabular}{|l|l|l|}
\hline
$(a,b)$&$D(f,g)$&\text{comment} \\
\hline
$(1,b)$ and $(a,1)$&$\int\otimes \id\colon  f'g$ and $fg'$&\text{for $\un=1$}\\
$(a,b)$&$P_{00}\colon \alpha f'g+\beta fg'$&\text{for any $\alpha,\beta\in\Kee$ and $\un=1$}\\
\hline
$(0,0)$&$P_{00}\colon \alpha f'g+\beta fg'$&\text{for any $\alpha,\beta\in\Kee$}\\
\hline
$(a,b)$&$\{\cdot,\cdot\}_{P.B.}\colon =afg'+bf'g$ &\text{Poisson bracket}\\
$(1,1)$&$\{\cdot,\cdot\}\colon =fg'+f'g$ &\text{contact bracket of generating functions}\\
\hline\end{tabular}
\end{equation}

$\bullet$ \emph{Order 2 operators} are the same as for $p=0$, but
the condition $a+b+1=0$ is to be replaced with $a+b+1=0 \pmod p$ and
for $p=2$ we have
\begin{equation}
\label{t2} \begin{tabular}{|l|l|l|}
\hline$(a,b)$&$D(f,g)$&\text{comment} \\
\hline
$(a,b)$ &$f'g'$&for $\un=1$; if $(a,b)=(1,1)$, then $D=\int\otimes \int$ \\
\hline
$(0,b)$ &$f'g'+bf''g$&$\{df,g(du)^b\}_{P.B.}$\\
$(a,1+a)$&$afg''+f'g'+(a+1)f''g$&$\{f(du)^a,dg(du)^{1+a}\}_{P.B.}$\\
$(a,b)$&$a g''f+f'g'+ b f'' g$&generalizes the above\\
\hline\end{tabular}
\end{equation}

whereas for $p=3$ we have
\begin{equation}
\label{t23} \begin{tabular}{|l|l|l|}
\hline$(a,b)$&$D(f,g)$&\text{comment} \\
\hline
$(a,b)$ for $a,b\not =0,1$&$J_{a,b}^2(f,g)$
&for $\un=1$\\
$(1,b)$ &$\int\otimes \id(f,g)= f''g$ &for $\un=1$\\
$(0,b)$ &$f'g' -b fg''$ &for $\un=1$\\
\hline

$(a,-1-a)$&$\{f,g\}'=afg''+(1-a)f'g'+(a+1)f''g$& \\
$(1,1)$&$Bj_{3,1}(f,g)$&\\
$(a,0)$&$\{f,g'\}=-afg''+f'g'$& \\
\hline\end{tabular}
\end{equation}

$\bullet$ \emph{Order 3 operators} are the same as for $p=0$, but
the Grozman operator $Gz$ does not survive for $p=2$ and $3$. For
$p=2$, we additionally have (this might be unclear; in fact, for
$(a,b)=(1,1)$, we have $\int\otimes \id$ and $\id \otimes \int$ and
$Bj_{2,2}=\int\otimes \id-\id \otimes \int$):
\begin{equation}
\label{t32} \renewcommand{\arraystretch}{1.4}
\begin{tabular}{|l|l|l|}
\hline
$(a,b)$&$D(f,g)$&\text{comment} \\
\hline
$(1,b)$ for $b\not =0$&$\int\otimes \id(f,g)= f'''g$&\text{for $\un=2$}\\
$(1,0)$&$Bj^{*2}_{2,2}$, $\int\otimes \id$ &\text{for $\un=2$}\\
$(a,b)$ for $a,b\not=0,1$&$J_{a,b}^3(f,g)$ 
&\text{for $\un=2$}\\
\hline
\end{tabular}
\end{equation}

$\bullet$ \emph{Order 4 operators}: no operators if $p>5$; for
$p\leq5$, we have

\begin{equation}
\label{t4} \begin{tabular}{|l|l|l|l|} \hline
$(a,b)$&$D(f,g)$&$p$&comments\\
\hline

$(0,0)$&$Gz_{2,2}(df,dg)=f'''g'+f'g''' +f'' g''$&$2$& $\un=2$\\
$(1,0)$&$Gz_{2,2}^{*1}(df,dg)=fg^{(4)} +f''' g'$&$2$& $\un=2$\\
$(0,1)$&$Gz_{2,2}^{*2}(df,dg)=f^{(4)} g+ f' g'''$&$2$& $\un=2$\\
\hline \hline

$(1,1)$& $\int\otimes\int=f''g''$&$3$&for $\un=1$\\
\hline
$(0,0)$& $Bj_{3,1}\circ (d f,dg)$ &$3$&\\
\hline
\hline

$(1,1)$& $Bj_{5,1}(f,g)$&$5$&\\
$(1,0)$&$Bj_{5,1}^{*2}(f,g)=
fg^{(4)}-f'''g'-f'g'''+f''g''$&$5$&\\
$(0,1)$&$Bj_{5,1}^{*1}(f,g)= f^{(4)}g-f'''g'+f''g''-f'g'''$&$5$&\\
\hline

$(1, b)$, $b\not = 0$&$\int \otimes \id \colon  f^{(4)}g$&$5$& for $\un=1$\\
$(0, 1)$  &$\id \otimes \int, Bj_{5,1}^{*1}$&$5$& for $\un=1$\\
$(a, 1-a), a\not =0,1$&$L_{5,1}(f,g)=f^{(4)} g+f g^{(4)}  -f''' g'
-f'g'''+f''g''$&$5$&\\
\hline\end{tabular}
\end{equation}

$\bullet$ \emph{Order 5 operators}: no invariant operators if $p>7$;
for $p\leq7$,  we have
\begin{equation}
\label{t5} \begin{tabular}{|l|l|l|l|} \hline
$(a,b)$&$D(f,g)$&$p$&comments\\
\hline
$(0,0)$& $Bj_{2,2}(df,dg)=f'g^{(4)}+f^{(4)}g'$ &$2$&\\
\hline

$(0,0)$& $Gz_{5,1}(df,dg)=f^{(4)} g'-g^{(4)}f'+f''g'''-f''' g''$ &$5$&\\
$(0,1)$& $Gz_{5,1}^{*2}(df,dg)=f^{(5)}g-f'g^{(4)}$&$5$&\\
$(1,0)$& $Gz_{5,1}^{*1}(df,dg)=fg^{(5)}-f^{(4)}g'$&$5$&\\
\hline

$(1,1)$&$Gz_{7,1}(f,g)$&$7$&\\
\hline\end{tabular}
\end{equation}

$\bullet$ \emph{Order 6 operators}: no invariant operators if $p>7$;
for $p\leq7$, we have
\begin{equation}
\label{t6} \begin{tabular}{|l|l|l|l|} \hline
$(a,b)$&$D(f,g)$&$p$&comment\\
\hline
$(1,1)$& $\int\otimes\int$ & $2$&for $\un=2$\\
$(1,1)$& $Gz_{2,3}(f,g)$ & $2$&for $\un>2$\\
 \hline
$(0,0)$& $Bj_{5,1}(df,dg)=f^{(5)} g'-f' g^{(5)}$& $5$&\\
\hline
$(0,1)$& $Bj_{7,1}^{*1}(f,g)$
& $7$&\\
$(1,1)$& $Bj_{7,1}(f,g) $& $7$&\\
$(1,0)$& $Bj_{7,1}^{*2}(f,g)$
& $7$&\\
\hline
$(1,b)$ for $b\not =0$& $\int\otimes\id$& $7$&for $\un=1$\\
$(1,0)$& $\int\otimes\id$, $Bj_{7,1}^{*2}$& $7$&for $\un=1$\\
$(0,1)$& $Bj_{7,1}^{*1}$, $\id\otimes\int$& $7$&for $\un=1$\\
$(1,1)$ & $Bj_{7,1}(f,g)$ & $7$&for $\un=1$\\
$(a,1-a)$, $a\not =0,1$& $L_{7,1}(f,g)$
& $7$&\\

\hline
\end{tabular}
\end{equation}

$\bullet$ \emph{Order 7 operators}: no invariant operators if $p>7$;
for $p\leq7$, we have
\begin{equation}
\label{t7} \begin{tabular}{|l|l|l|l|} \hline
$(a,b)$&$D(f,g)$&$p$&comment\\
\hline $(1,0)$ & $Bj_{2,3}^{*2}(f,g)$ 
&$2$&\\

$(1,1)$ &$Bj_{2,3}(f,g)$ &$2$&\\

$(1,b)$ for $b\not = 0$& $\int \otimes \id $& $2$& $N=3$\\
$(1,0)$& $\int \otimes \id$, $Bj_{2,3}^{*2}$& $2$& $N=3$\\
$(a,1-a), a\not =0,1$ &$L_{2,3}(f,g)$
& $2$& $N=3$\\

\hline $(1,1)$ & $Gz_{3,2}(f,g)$ &$3$&\\
\hline

$(0,0)$ &$Gz_{7,1}(df,dg)$
& $7$ &\\
$(1,0)$ &$Gz_{7,1}^{*1}(df,dg)=fg^{(7)}-f^{(6)} g'$& $7$ &\\
$(0,1)$ &$Gz_{7,1}^{*2}(df,dg)=f^{(7)} g-f'g^{(6)}$& $7$ &\\
\hline
\end{tabular}
\end{equation}

\subsection*{Acknowledgements} We are thankful to Grozman for his wonderful
code \textit{SuperLie}, see \cite{Gr}. S.B. was partly supported by
the grant AD 065 NYUAD.

\end{document}